\documentclass[12pt,twoside]{article}
\usepackage[margin=0.8in]{geometry}
\usepackage{amsmath, amssymb, amsfonts, mathtools}
\usepackage{amsthm}
\usepackage{mathrsfs}
\usepackage{stmaryrd}
\usepackage{esint}
\usepackage{eqnalign}
\usepackage{graphicx}
\usepackage{enumitem}
\usepackage[most]{tcolorbox}
\usepackage{hyperref}
\usepackage{float}
\usepackage{listings}
\usepackage{array}
\usepackage{booktabs}
\usepackage{xcolor}
\usepackage{emptypage}
\usepackage{tensor}
\usepackage{array}
\usepackage{tabularx}
\usepackage{hyperref}
\usepackage[nameinlink,noabbrev]{cleveref}
\usepackage{titlesec}
\usepackage{fancyhdr}
\usepackage[T1]{fontenc}
\usepackage[utf8]{inputenc}
\usepackage{lmodern}

\usepackage{titling}   
\usepackage{fancyhdr}
\usepackage[backend=biber,style=numeric]{biblatex}

\newcommand{\EulerB}[2]{\left\langle {#1 \atop #2} \right\rangle^{\!B}}

\newcommand{\HeadTitle}{Perspectives on the arithmetic nature of the ratios 
  \ensuremath{\zeta(2n+1)/\pi^{2n+1}} and 
  \ensuremath{\beta(2n)/\pi^{2n}}}
\newcommand{\HeadAuthor}{Luc Ramsès TALLA WAFFO}

\titleformat{\section} {\normalfont\normalsize\bfseries\itshape\centering} {\thesection} {0em} {}

\titleformat{\section}[block]
  {\normalfont\large\bfseries\itshape\centering}
  {§\thesection.}
  {1em}
  {}

\titleformat{\subsection}
  {\normalfont\normalsize\bfseries}
  {\thesubsection}
  {1em}
  {}
	\newtheorem{theorem}{Theorem}[section]

\newtheorem{proposition}[theorem]{Proposition}

\crefname{lemma}{lemma}{lemmas}
\Crefname{lemma}{Lemma}{Lemmas}

\crefname{proposition}{proposition}{propositions}
\Crefname{proposition}{Proposition}{Propositions}

\crefname{theorem}{theorem}{theorems}
\Crefname{theorem}{Theorem}{Theorems}

\begin{document}

\thispagestyle{fancy} 

\begin{center}
\textbf{}
\end{center}

\begin{center}
\Large{\textit{Perspectives on the arithmetic nature of the ratios}}\end{center} 
\begin{center}
\Large{\textit{$\dfrac{\zeta(2n+1)}{\pi^{2n+1}}$ and $\dfrac{\beta(2n)}{\pi^{2n}}$}}
\end{center}

\hspace{3cm}
\begin{center}
Luc Ramsès TALLA WAFFO \\
Technische Universität Darmstadt\\
Karolinenplatz 5, 64289 Darmstadt, Germany\\
ramses.talla@stud.tu-darmstadt.de\\
\vspace{0.5cm}
October 8, 2025
\end{center}

\begin{abstract} We investigate special values of the Riemann zeta function at odd integers and of the Dirichlet beta function at even integers. We collect several analytic frameworks leading to these values and thereby provide a unifying perspective. Beyond their analytic interest, the resulting formulae motivate linear-independence conjectures which, if proved, would imply the irrationality of $\zeta(2n+1)/\pi^{2n+1}$ and $\beta(2n)/\pi^{2n}$.
\end{abstract}

\vspace{0.2cm}

\section*{Introduction}

\vspace{0.5cm}

Euler is renowned for his resolution of the Basel problem, establishing that  
\[
\zeta(2) = \frac{\pi^2}{6}.
\]  
Moreover, he derived a general formula for the values of the Riemann zeta function at even positive integers, namely  
\[
\zeta(2n) = (-1)^{\,n+1}\, \frac{B_{2n}\,(2\pi)^{2n}}{2\,(2n)!},
\]  
where \(B_{2n}\) denotes the \(2n\)-th Bernoulli number \cite{Euler1734}. Despite his extensive efforts, Euler was unable to obtain an analogous formula for \(\zeta(2n+1)\). Since then, the arithmetic nature of these odd zeta values has remained a central open problem, attracting sustained attention from mathematicians to this day.  

\vspace{0.1cm}

A closely related, though less widely studied, function is \emph{Dirichlet’s beta function}, defined for \(\Re(s) > 0\) by  
\[
\beta(s) := \sum_{m=0}^\infty \frac{(-1)^m}{(2m+1)^s}.
\]  
Although this \(L\)-function shares many structural similarities with \(\zeta(s)\), it has received comparatively limited attention in the literature. Historically, even prior to Euler’s investigations, Leibniz demonstrated the identity  
\[
\sum_{n=0}^\infty \frac{(-1)^n}{2n+1} \;=\; \frac{\pi}{4}.
\]  

In contrast to the Riemann zeta function, whose explicit values are known at even arguments but remain mysterious at odd ones, the Dirichlet beta function exhibits a complementary phenomenon. For odd arguments, one has explicit closed‐form evaluations \cite{Zudilin2019}:  
\[
\beta(2n+1) = (-1)^n \,\frac{E_{2n}\,\pi^{2n+1}}{2^{2n+2}\,(2n)!},
\]  
where \(E_{2n}\) are Euler numbers. Hence, the ratios \(\beta(2n+1)/\pi^{2n+1}\) are rational, underscoring a striking parallel with the special values of \(\zeta(s)\). For even arguments of \(\beta(s)\), however, no such closed formula is currently known. 

\vspace{0.1cm}

Motivated by Euler's evaluation of \(\zeta(2n)\) and the expression of $\beta(2n+1)$, many authors tried to find some rational coefficients $r_n$ and $s_n$ such that the given relations may be satisfied
\[\zeta(2n + 1) = r_n \, \pi^{2n+1} \quad \text{and } \quad \beta(2n) = s_n \, \pi^{2n}\]

Given the futility of the extraordinarily numerous attempts to compute these numbers, the solution of the designated problem is generally regarded as impossible; but so far a strong proof of this impossibility is lacking. Only heuristic methods and deep computations by brute force support this impossibility to this day.  \textit{One may solve the problem by proving that the numbers $\dfrac{\zeta(2n+1)}{\pi^{2n}}$ and $\dfrac{\beta(2n)}{\pi^{2n-1}}$ are not multiples of $\pi$.} The difficulty of this new reformulation partly stems from the lack
of sufficiently many closed-form expressions or identities involving these values that might
serve as a starting point for an irrationality or transcendence proof. \textit{This paper is devoted to some integral representations of the numbers $\dfrac{\zeta(2n+1)}{\pi^{2n}}$ and $\dfrac{\beta(2n)}{\pi^{2n-1}}$ through interesting logarithmic, trigonometric and hyperbolic kernels, along with several conjectures that support their non-representability as rational multiples of $\pi$}.

\vspace{0.5cm}

\section{\hspace{0.2cm}Preliminaries}\label{sec:preliminaries}

\vspace{0.3cm}

The point of departure of our investigations was the Malmsten's integral \[I = \displaystyle\int_0^1\dfrac{x(x^4-4x^2+1)\ln\ln\frac{1}{x}}{(1+x^2)^4}dx\] that was first evaluated by \textit{Iaroslav V, Blagouchine} in 2014 \cite{Blagouchine2014} and shown to equal $\dfrac{7\zeta(3)}{8\pi^2}$. Therefore the proof of its convergence is omitted here for the sake of conciseness. Notwithstanding, the presence of  singularities at its bounds requires some care. Consequently, we rewrite the integral as the limit $I = \displaystyle\lim_{\xi \to 0^+} \int_\xi^{1-\xi}\dfrac{x(x^4-4x^2+1)\ln\ln\frac{1}{x}}{(1+x^2)^4}dx$

\vspace{0.3cm}

We consider first $u'(x) = \dfrac{x(x^4-4x^2+1)}{(1+x^2)^4}$ and we show easily that its anti-derivative is $u(x) = \dfrac{x^2(1-x^2)}{2(1+x^2)^3}$. Then we consider $v(x) = \ln\ln\frac{1}{x}$ and observe that its derivative is $v'(x) = \dfrac{1}{x \ln x}$. Since both functions $u, v \in \mathcal C^1[\xi, 1 - \xi]$, an integration by parts yields

\[ I = \displaystyle\lim_{\xi \to 0^+} \left(\left[\dfrac{x^2(1-x^2)\ln\ln\frac{1}{x}}{2(1+x^2)^3}\right]_\xi^{1-\xi}-\dfrac{1}{2}\int_\xi^{1-\xi}\dfrac{x^2(1-x^2)}{(1+x^2)^3x\ln x}dx\right) \]

This implies 
\[ I =  \displaystyle\underbrace{\lim_{\xi \to 0^+} \left[\dfrac{x^2(1-x^2)\ln\ln\frac{1}{x}}{2(1+x^2)^3}\right]_\xi^{1-\xi}}_{0}-\lim_{\xi \to 0^+}\dfrac{1}{2}\int_\xi^{1-\xi}\dfrac{x^2(1-x^2)}{(1+x^2)^3x\ln x}dx \]

It remains 
\[ I = \displaystyle\lim_{\xi \to 0^+}\dfrac{1}{2}\int_\xi^{1-\xi}\dfrac{x(x^2-1)}{(1+x^2)^3\ln x}dx = \dfrac{1}{2}\int_0^1\dfrac{x(x^2-1)}{(1+x^2)^3\ln x}\,dx \]

The change of variable $x \rightarrow \tan u$ yields:
\begin{allowdisplaybreaks}
\begin{align*}
\displaystyle\int_0^1\dfrac{x(x^2-1)}{(1+x^2)^3\ln x}\,dx = \dfrac{7\zeta(3)}{4\pi^2} & \Rightarrow \displaystyle\int_0^{\frac{\pi}{4}}\dfrac{\tan u(\tan^2u-1)}{(1+\tan^2u)^3\ln( \tan u)}(1+\tan^2u)\,du = \dfrac{7\zeta(3)}{4\pi^2} \\
& \Rightarrow \displaystyle\int_0^{\frac{\pi}{4}}\dfrac{\cos^4u\tan u(\tan^2u-1)}{\ln( \tan u)}\,du = \dfrac{7\zeta(3)}{4\pi^2}\\
& \Rightarrow \displaystyle\int_0^{\frac{\pi}{4}}\dfrac{\cos u\sin u(\sin^2u-\cos^2u)}{\ln( \tan u)}\,du = \dfrac{7\zeta(3)}{4\pi^2} \\
& \Rightarrow -\displaystyle\int_0^{\frac{\pi}{4}}\dfrac{\cos u\sin u \cos(2u)}{\ln( \tan u)}\,du = \dfrac{7\zeta(3)}{4\pi^2}  
\end{align*}
\end{allowdisplaybreaks}

This gives finally 
 
\begin{equation}
\displaystyle\int_0^{\frac{\pi}{4}}\dfrac{\sin 4u }{\ln( \tan u)}\,du = -\dfrac{7\zeta(3)}{\pi^2}\label{eq:zeta3}
\end{equation}

\vspace{0.3cm}

Blagouchine introduced in fact more generally the family of integrals on page 80:

\[I_n = \int_0^1 \frac{x^{n-1} \ln \ln \frac{1}{x}}{(1 + x^2)^n} \, dx 
= \int_1^\infty \frac{x^{n-1} \ln \ln x}{(1 + x^2)^n} \, dx 
= \frac{1}{2n} \int_0^\infty \frac{\ln x}{\cosh^n x} \, dx,\]

And with contour integration technique, he obtained further\cite{Blagouchine2014} :

\[
I_2 = -\frac{1}{2} \ln 2 + \frac{1}{4} \ln \pi - \frac{\gamma}{4}
\]

\[
I_6 = -\frac{1}{60} \ln 2 + \frac{1}{120} \ln \pi - \frac{\gamma}{120}
      - \frac{7\,\zeta(3)}{192\pi^2} - \frac{31\,\zeta(5)}{320\pi^4}
\]

where $\gamma$ denotes the Euler-Mascheroni constant. We can rewrite $I_6$ as \[
I_6 = -\frac{1}{30} I_2 - \frac{7\,\zeta(3)}{192\pi^2} - \frac{31\,\zeta(5)}{320\pi^4}
\]

And by replacing $I_6, I_2$ and $\zeta(3)$ by their $\ln \ln$ integral representations, it follows:

\[\int_0^1 \dfrac{x(x^8-26x^6+66x^4-26x^2+1)\ln \ln \frac{1}{x}}{(1+x^2)^6}\,dx = -\dfrac{93\zeta(5)}{8\pi^4}\]
Same remarks and same algebra on $I_8$ yield:

\[\int_0^1 \dfrac{x(x^{12}-120x^{10}+1\,191x^8-2\,416x^6+1\,191x^4-120x^2+1)\ln \ln \frac{1}{x}}{(1+x^2)^8}\,dx = \dfrac{5\,715\zeta(7)}{16\pi^6}\]

And by pursuing the same line of reasoning as above, the two later Malmsten's integrals take the forms :
\begin{equation}\dfrac{\zeta(5)}{\pi^{4}} = -\dfrac{1}{186}\int_0^{\frac{\pi}{4}}\dfrac{\sin 4x}{\ln(\tan x)}\,dx +\dfrac{1}{124}\int_0^{\frac{\pi}{4}}\dfrac{\sin 8x}{\ln(\tan x)}\,dx \label{eq:zeta5}\end{equation}
\begin{equation}
\frac{\zeta(7)}{\pi^{6}} = -\frac{17}{91\,440} \int_0^{\frac{\pi}{4}} \frac{\sin(4x)}{\ln(\tan x)} \, dx
+ \frac{1}{1\,524} \int_0^{\frac{\pi}{4}} \frac{\sin(8x)}{\ln(\tan x)} \, dx
- \frac{1}{2\,032} \int_0^{\frac{\pi}{4}} \frac{\sin(12x)}{\ln(\tan x)} \, dx \label{eq:zeta7}\end{equation}
\vspace{0.5cm}

These results highlight general patterns. In fact, it is easy to recognize the pattern \[\displaystyle\int_0^{\frac{\pi}{4}} \dfrac{\sin(4nx)}{\ln(\tan x)} \, dx\] where $n \in \mathbb{N}$. 

\vspace{0.5cm}

\section{\hspace{0.2cm}Main Result}\label{sec:result}

\vspace{0.5cm}

Deeper linear algebra through the Gauss elimination algorithm on \eqref{eq:zeta3}, \eqref{eq:zeta5} and \eqref{eq:zeta7} allows us to conclude the following :

\[
\int_0^{\frac{\pi}{4}} \dfrac{\sin(4x)}{\ln(\tan x)} \, dx = -7 \dfrac{\zeta(3)}{\pi^2}
\]

\[
\int_0^{\frac{\pi}{4}} \dfrac{\sin(8x)}{\ln(\tan x)} \, dx = -\dfrac{14}{3} \dfrac{\zeta(3)}{\pi^2} + 124 \dfrac{\zeta(5)}{\pi^4}
\]

\[
\int_0^{\frac{\pi}{4}} \dfrac{\sin(12x)}{\ln(\tan x)} \, dx = -\dfrac{161}{45} \dfrac{\zeta(3)}{\pi^2} + \dfrac{496}{3} \dfrac{\zeta(5)}{\pi^4} - 2\,032 \dfrac{\zeta(7)}{\pi^6}
\]

After further computations, the following equalities also hold true :

\[
\int_0^{\frac{\pi}{4}} \dfrac{\sin 16x}{\ln(\tan x)} \, dx = -\dfrac{44}{15} \dfrac{\zeta(3)}{\pi^2}
+ \dfrac{2\,728}{15} \dfrac{\zeta(5)}{\pi^4}
- 4\,064 \dfrac{\zeta(7)}{\pi^6}
+ 32\,704 \dfrac{\zeta(9)}{\pi^8}
\]

\[
\int_0^{\frac{\pi}{4}} \dfrac{\sin 20x}{\ln(\tan x)} \, dx =
-\dfrac{563}{225} \dfrac{\zeta(3)}{\pi^2}
+ \dfrac{178\,064}{945} \dfrac{\zeta(5)}{\pi^4}
- \dfrac{87\,376}{15} \dfrac{\zeta(7)}{\pi^6}
+ \dfrac{261\,632}{3} \dfrac{\zeta(9)}{\pi^8}
- 524\,032 \dfrac{\zeta(11)}{\pi^{10}}
\]

\begin{align*}
\int_0^{\frac{\pi}{4}} \dfrac{\sin 24x}{\ln(\tan x)} \, dx &=
-\dfrac{3\,254}{1\,485} \dfrac{\zeta(3)}{\pi^2}
+ \dfrac{2\,697\,868}{14\,175} \dfrac{\zeta(5)}{\pi^4}
- \dfrac{1\,381\,760}{189} \dfrac{\zeta(7)}{\pi^6} 
+ \dfrac{6\,933\,248}{45} \dfrac{\zeta(9)}{\pi^8}
- \dfrac{5\,240\,320}{3} \dfrac{\zeta(11)}{\pi^{10}} \\
&+ 8\,387\,584 \dfrac{\zeta(13)}{\pi^{12}}
\end{align*}

\begin{align*}
\int_0^{\frac{\pi}{4}} \dfrac{\sin 28x}{\ln(\tan x)} \, dx =\;&
-\dfrac{88\,069}{45\,045} \dfrac{\zeta(3)}{\pi^2}
+ \dfrac{1\,409\,632}{7\,425} \dfrac{\zeta(5)}{\pi^4}
- \dfrac{8\,091\,424}{945} \dfrac{\zeta(7)}{\pi^6} \\
&+ \dfrac{30\,685\,696}{135} \dfrac{\zeta(9)}{\pi^8}
- 3\,668\,224 \dfrac{\zeta(11)}{\pi^{10}}
+ 33\,550\,336 \dfrac{\zeta(13)}{\pi^{12}} \\
&- 134\,213\,632 \dfrac{\zeta(15)}{\pi^{14}}
\end{align*}

\begin{align*}
\int_0^{\frac{\pi}{4}} \dfrac{\sin 32x}{\ln(\tan x)} \, dx =\;&
-\dfrac{11\,384}{6\,435} \dfrac{\zeta(3)}{\pi^2}
+ \dfrac{4\,448\,479\,664}{23\,648\,625} \dfrac{\zeta(5)}{\pi^4}
- \dfrac{300\,061\,376}{31\,185} \dfrac{\zeta(7)}{\pi^6} \\
&+ \dfrac{614\,115\,712}{2\,025} \dfrac{\zeta(9)}{\pi^8}
- \dfrac{391\,975\,936}{63} \dfrac{\zeta(11)}{\pi^{10}}
+ \dfrac{1\,224\,587\,264}{15} \dfrac{\zeta(13)}{\pi^{12}} \\
&- \dfrac{1\,878\,990\,848}{3} \dfrac{\zeta(15)}{\pi^{14}}
+ 2\,147\,467\,264 \dfrac{\zeta(17)}{\pi^{16}}
\end{align*}

These considerations motivate the following conjecture:
 \[\forall n \in \mathbb{N}^*, \exists \, \begin{pmatrix}
C_{1,n} \\
C_{2,n} \\
\vdots  \\
C_{n,n} 
\end{pmatrix} \in \mathbb{Q}^{n} \quad \text{such that} \quad \int_0^{\frac{\pi}{4}} \dfrac{\sin(4nx)}{\ln(\tan x)} \, dx  = \sum_{p=1}^{n} C_{p,n} \, \frac{\zeta(2p+1)}{\pi^{2p}}\]

The proof of this statement is fully detailed in \cite[59-67]{talla_waffo_integral_2025} and is omitted here in the interest of brevity. We merely note the existence of integrals $I_n$ satisfying
\begin{equation}\forall n \in \mathbb{N}^*, \exists \, \begin{pmatrix}
C_{1,n} \\
C_{2,n} \\
\vdots  \\
C_{n,n} 
\end{pmatrix} \in \mathbb{Q}^{n} \quad \text{such that} \quad I_n  = \sum_{p=1}^{n} C_{p,n} \, \frac{\zeta(2p+1)}{\pi^{2p}} \label{eq:formZeta}\end{equation}

Furthermore, we also find in \cite[49-86]{talla_waffo_integral_2025} some integrals $I_n$ which adhere to this form

\begin{equation}\forall n \in \mathbb{N}^*, \exists \, \begin{pmatrix}
C_{1,n} \\
C_{2,n} \\
\vdots  \\
C_{n,n} 
\end{pmatrix} \in \mathbb{Q}^{n} \quad \text{such that} \quad I_n  = \sum_{p=1}^{n} C_{p,n} \, \frac{\beta(2p)}{\pi^{2p-1}} \label{eq:formBeta}\end{equation}

\textit{The formal existence of integrals of the type \eqref{eq:formZeta} or \eqref{eq:formBeta} can already be traced back to the work of Tobias Kyrion \textnormal{\cite{Kyrion2025}}, where the author investigates, among others, the following integrals:}
\[
\int_{0}^{+\infty} \frac{\tanh x}{x}\, \operatorname{sech}^{2n+1}\! x\, dx, 
\qquad 
\int_{0}^{+\infty} \frac{\tanh x}{x}\, \operatorname{sech}^{2n}\! x\, dx,
\qquad 
\int_{0}^{+\infty} \frac{\tanh^{n} x}{x^{n}}\, dx.
\]
\textit{These three families of integrals were later rediscovered in \textnormal{\cite{talla_waffo_integral_2025}} within a broader analytic framework, and are adopted here in a different formulation that is more suitable for the evaluation method and for discussing certain linear independence conjectures. Moreover, several additional integrals also satisfy either \eqref{eq:formZeta} or \eqref{eq:formBeta}. The collection of all such integrals, together with the remarkable structural patterns they exhibit, constitutes the main contribution of this article. A representative selection is summarized in the following table, where $n$ denotes a strictly positive integer and $C_{p,n}$ are rational coefficients with $p \in \llbracket 1, n \rrbracket$.}

\vspace{0.5cm}

\begin{tabularx}{\textwidth}{|>{\centering\arraybackslash}X|
                               >{\centering\arraybackslash}X|}
\hline
$\displaystyle \sum_{p=1}^{n} C_{p,n} \, \frac{\zeta(2p+1)}{\pi^{2p}}$ 
& 
$\displaystyle \sum_{p=1}^{n} C_{p,n} \, \frac{\beta(2p)}{\pi^{2p-1}}$ \\
\hline
\begin{itemize}
  \item $\tiny{\displaystyle \int_{0}^{+\infty} 
          \frac{\sinh\!\left((2k+1)x\right)}{x \cosh^{2n+1}\! x} \, dx 
          \quad (k \in \llbracket 0, n-1 \rrbracket)}$
  \item $\displaystyle \int_{0}^{+\infty} 
          \frac{\sinh(2k x)}{x \cosh^{2n+2} x} \, dx 
          \quad (k \in \llbracket 1, n \rrbracket)$
  \item $\displaystyle \int_{0}^{\tfrac{\pi}{4}} 
          \frac{\sin(4n x)}{\ln(\tan x)} \, dx$
  \item $\displaystyle \int_{0}^{1} 
          \frac{\operatorname{Li}_{-2n-1}(-x^{2}) }{x}\, \ln \ln \tfrac{1}{x} \, dx$
  \item $\displaystyle \int_{0}^{1} 
          \frac{x^{2n-1}}{\operatorname{arctanh} x} \, dx$
  \item $\displaystyle \int_{0}^{+\infty} 
          \frac{\tanh^{2n} x}{x^{2}} \, dx$
  \item $\displaystyle \int_{0}^{+\infty} 
          \frac{\tanh^{\,n+1} x}{x^{\,n+1}} \, dx$
					\vspace{0.5cm}
\end{itemize}
&
\begin{itemize}
  \item $\tiny{\displaystyle \int_{0}^{+\infty} 
          \frac{\sinh\!\left((2k+1)x\right)}{x \cosh^{2n}x} \, dx 
          \quad (k \in \llbracket 0, n-1 \rrbracket)}$
  \item $\displaystyle \int_{0}^{+\infty} 
          \frac{\sinh(2k x)}{x \cosh^{2n+1} x} \, dx 
          \quad (k \in \llbracket 1, n \rrbracket)$
  \item $\displaystyle \int_{0}^{\tfrac{\pi}{4}} 
          \frac{\cos\!\big((4n-2)x\big)}{\ln(\tan x)} \, dx$
  \item $\displaystyle \int_{0}^{1} 
          \frac{\operatorname{Im}\!\left[\operatorname{Li}_{-2n}(i x)\right]}{x} 
          \, \ln \ln \tfrac{1}{x} \, dx$
  \item $\displaystyle \int_{0}^{1} 
          \frac{x^{2n-1}}{\sqrt{1-x^{2}} \, \operatorname{arctanh} x} \, dx$
\end{itemize} \\
\hline
\end{tabularx}

\vspace{0.5cm}

The reader is referred to \textnormal{\cite{talla_waffo_integral_2025}} for illustrative examples of these integrals and for detailed derivations. We remark in passing that the integral $\displaystyle \int_{0}^{1} 
          \frac{\operatorname{Im}\!\left[\operatorname{Li}_{-2n}(i x)\right]}{x} 
          \, \ln \ln \tfrac{1}{x} \, dx$ generalizes the integral \[
\dfrac{\beta(2)}{\pi}
= \displaystyle\int_0^1 \dfrac{u^4 - 6u^2 + 1}{2(1 + u^2)^3} \,\ln \ln \frac{1}{u} \,du
\] that has been firstly evaluated by Adamchik \cite{Adamchik1997} and rediscovered by Blagouchine \cite{Blagouchine2014}. Additionally, the integral $\displaystyle \int_{0}^{1} \frac{\operatorname{Li}_{-2n-1}(-x^{2}) }{x}\, \ln \ln \tfrac{1}{x} \, dx$ exhibits a general form of $\dfrac{\zeta(2n+1)}{\pi^{2n}}$ by also extending the Blagouchine's integral $\displaystyle\int_0^1\dfrac{x(x^4-4x^2+1)\ln\ln\frac{1}{x}}{(1+x^2)^4}dx = \frac{7\, \zeta(3)}{8\, \pi^2}$ to all higher $n$. Strictly speaking, the endeavour to generalize the integral representations of Adamchik and Blagouchine did not originate with the present work. Indeed, Tobias Kyrion, in his article \textnormal{\cite{Kyrion2025}}, establishes the following formulae, where $N$ denotes a strictly positive integer: 

\[
\begin{aligned}
\beta(2N)
&= (-1)^{N}\,
\frac{\pi^{\,2N-1}}{2^{\,2N-1}(2N-1)!}
\int_{1}^{\infty}
\sum_{k=0}^{N-1}
\left(
  \sum_{j=0}^{k}
  \binom{2k+1}{k-j}
  (-1)^{\,j+1}
  (2j+1)^{\,2N-1}
\right)
\\[1em]
&\quad \times
\frac{
  \bigl((2k+1)x^{4}-2(2k+3)x^{2}+2k+1\bigr)\,x^{2k}
}{
  (x^{2}+1)^{2k+3}
}
\log(\log(x))\, dx .
\end{aligned}
\]

\hspace{0.3cm}

\[
\begin{aligned}
\zeta(2N+1)
&= (-1)^{N}\,
\frac{2\pi^{\,2N}}{(2^{\,2N+1}-1)\,(2N)!}
\int_{1}^{\infty}
\sum_{k=1}^{N}
\left(
  \sum_{j=1}^{k}
  \binom{2k}{k-j}
  (-1)^{j}
  (2j)^{\,2N}
\right)
\\[1em]
&\quad \times
\frac{
  \bigl(2k x^{4} - 2(2k+2)x^{2} + 2k\bigr)\, x^{2k-1}
}{
  (x^{2}+1)^{2k+2}
}
\log(\log(x))\, dx .
\end{aligned}
\]

Although these formulae are elegant and reveal the underlying structure of the integral, they involve multiple nested summations together with binomial coefficients, which makes them rather cumbersome to handle in irrationality proofs.
The general integrals listed in the table above, by contrast, yield new representations that avoid such complications entirely, as the only special function that appears is the well-known Jonquière's polylogarithm.
Closed-form expressions for these new integrals are provided in the following propositions.

\vspace{0.3cm}

\begin{proposition}\label[proposition]{prop:zeta_polylog}
\[
\forall \, n \in \mathbb{N}^*, \quad
\int_0^1 \frac{\mathrm{Li}_{-2n-1}(-x^2) \, \ln\ln \frac{1}{x}}{x} \, dx
= (-1)^n \left(1 - \frac{1}{2^{2n+1}}\right) \,
\dfrac{(2n)!}{2}\frac{\zeta(2n+1)}{\pi^{2n}} 
\]
\end{proposition}

\begin{proof}

Let $J_n := \displaystyle\int_0^1 \dfrac{\mathrm{Li}_{-2n-1}(-x^2) \ln \ln \frac{1}{x}}{x}\,dx$. Recalling the definition of the polylogarithm $\displaystyle\mathrm{Li}_s(z) := \sum_{k = 1}^{\infty}\dfrac{z^k}{k^s}$ and bearing in mind that $|-x^2| < 1$ on the domain of integration, one gets

\[J_n = \int_0^{1} \dfrac{1}{x}\,\displaystyle\sum_{k = 1}^{\infty}\dfrac{(-1)^k\,x^{2k}}{k^{-2n-1}} \ln \ln \frac{1}{x}\,dx = \int_0^1 \displaystyle\sum_{k = 1}^{\infty} (-1)^k\,x^{2k-1}k^{2n+1} \ln \ln \frac{1}{x}\,dx\]

After switching the order of the operators (under suitable regularization assumptions), one has

\[J_n  = \displaystyle\sum_{k = 1}^{\infty} (-1)^k k^{2n+1} \int_0^1 x^{2k-1}\ln \ln \frac{1}{x}\,dx\]

And applying the well-known formula $\displaystyle\int_0^1 x^s \ln \left(\ln \frac{1}{x}\right)\,dx =-\dfrac{\gamma + \ln (s+1)}{s+1}$  gives 

\[J_n  = -\displaystyle\sum_{k = 1}^{\infty} (-1)^k k^{2n+1} \dfrac{\gamma + \ln (2k)}{2k} = -\dfrac{1}{2}\displaystyle\sum_{k = 1}^{\infty} (-1)^k k^{2n+1} \dfrac{\ln\left(2e^{\gamma}\right) + \ln k}{k}\]

The latter sum is then split in two to provide

\[J_n  =  -\dfrac{\ln\left(2e^{\gamma}\right)}{2}\underbrace{\displaystyle\sum_{k = 1}^{\infty} (-1)^k \dfrac{1}{k^{-2n}}}_{\eta(-2n)} - \dfrac{1}{2} \underbrace{\displaystyle\sum_{k = 1}^{\infty} (-1)^k \dfrac{\ln k}{k^{-2n}}}_{\eta'(-2n)}\]

One recognizes easily the Dirichlet $\eta$ function and its derivative. Now, the $\eta$-function has zeros at negative even integers since the Riemann $\zeta$ function also has zeros at negative even integers \cite{Riemann1859}\cite{WikiZetaFunctionalEquation}\cite{Euler1768}. And its derivative is \[\eta'(s) = \left(1 - 2^{1 - s}\right)\zeta'(s) + 2^{1 - s} \ln 2 \cdot \zeta(s)\]. It remains

\[ J_n = -\dfrac{1}{2} \left(1 - 2^{1 + 2n}\right)\zeta'(-2n)\]

And the closed form of the derivative zeta at even negative integers \cite{WikiZetaDerivatives}\cite{Apostol1976} \[\zeta'(-2n)
=
(-1)^n \,\frac{(2n)!}{2\,(2\pi)^{2n}}\,\zeta(2n+1)\]  yields

\[ J_n = (-1)^n \left(1 - \frac{1}{2^{2n+1}}\right) \,
\dfrac{(2n)!}{2}\frac{\zeta(2n+1)}{\pi^{2n}}\] This ends the proof.
\end{proof}

\begin{proposition}\label[proposition]{prop:beta_polylog}
\[\forall n \in \mathbb{N}^*,\int_0^1 
\dfrac{\operatorname{Im}\!\bigl[\operatorname{Li}_{-2n}(i\,x)\bigr]}{x} \,\ln \ln \frac{1}{x} \,dx = (-1)^{n+1} \dfrac{2^{2n-1}(2n-1)!}{\pi^{2n-1}}\beta(2n)\]
\end{proposition}

\begin{proof}
Malmsten proved in 1842 \cite{Blagouchine2014} the following functional equation

\[\beta(1 - s) =
\left( \frac{\pi}{2} \right)^{-s}
\sin\left( \frac{\pi s}{2} \right)
\Gamma(s)\, \beta(s)\]

Let $\chi(s) := \left( \dfrac{\pi}{2} \right)^{-s}
\sin\left( \dfrac{\pi s}{2} \right)\Gamma(s)$, so we get $\beta(1 - s) = \chi(s) \, \beta(s)$. Differentiating with respect to $s$, we get
\[-\beta'(1-s) = \chi'(s) \, \beta(s) + \chi(s) \, \beta'(s)\]

We see trivially that $\chi(2n) = 0$ due to the $\sin$ function. It remains

\[-\beta'(1-2n) = \chi'(2n) \, \beta(2n)\]

One has \[\chi'(2n) = \left. \dfrac{d}{ds} \left(\chi(s)\right)\right|_{s=2n}\]. Let $\Delta(s) := \left( \dfrac{\pi}{2} \right)^{-s} \Gamma(s)$, so that $\chi(s) = \Delta(s) \sin\left( \dfrac{\pi s}{2} \right)$. Differentiating with respect to $s$, we have

\[\chi'(s) = \Delta'(s) \sin\left( \dfrac{\pi s}{2} \right) + \dfrac{\pi}{2} \Delta(s) \cos\left( \dfrac{\pi s}{2} \right)\]

Evaluating this equality for $s = 2n$ and putting it back in the earlier relation, we obtain 
\begin{equation}
\beta'(1-2n) = (-1)^{n+1} \dfrac{2^{2n-1}(2n-1)!}{\pi^{2n-1}}\beta(2n) \label{eq:beta_derivatives_odd}\end{equation}

Let $R_n := \displaystyle\int_0^1 
\dfrac{\operatorname{Im}\!\bigl[\operatorname{Li}_{-2n}(i\,x)\bigr]}{x} \,\ln \ln \frac{1}{x} \,dx$. One has $\dfrac{\operatorname{Im}\!\bigl[\operatorname{Li}_{-2n}(i\,x)\bigr]}{x} = \displaystyle\sum_{k=0}^{\infty} (-1)^k\,(2k+1)^{2n} x^{2k}$ 

The Malmsten's integral becomes this after interchanging summation and integration (under suitable convergence assumptions)

\[R_n = \int_0^1 
 \sum_{k=0}^{\infty} (-1)^k\,(2k+1)^{2n} x^{2k}\,\ln \ln \frac{1}{x} \,dx = \sum_{k=0}^{\infty} (-1)^k\,(2k+1)^{2n}\int_0^1 
  x^{2k}\,\ln \ln \frac{1}{x} \,dx\]
	
After evaluating the integral, we get
	
	\[\int_0^1 
  x^{2k}\,\ln \ln \frac{1}{x} \,dx = -\dfrac{\gamma + \ln(2k+1)}{2k+1}\]
	
	We obtain
	
	\[R_n = -\sum_{k=0}^{\infty} (-1)^k\,(2k+1)^{2n} \dfrac{\gamma + \ln(2k+1)}{2k+1}\]

This yields

\[R_n = -\gamma\underbrace{\sum_{k=0}^{\infty} (-1)^k\,(2k+1)^{2n-1}}_{\beta(1-2n)} - \underbrace{\sum_{k=0}^{\infty} (-1)^k\,(2k+1)^{2n}\dfrac{\ln(2k+1)}{2k+1}}_{-\beta'(1-2n)}\]

One reads easily from the functional equation of $\beta$ that $\beta(1-2n) = 0$. The value of $\beta'(1-2n)$ follows from \eqref{eq:beta_derivatives_odd}.  

\vspace{0.3cm}

We conclude

\[R_n = \beta'(1-2n) = (-1)^{n+1} \dfrac{2^{2n-1}(2n-1)!}{\pi^{2n-1}}\beta(2n)\]This concludes the proof.
\end{proof}

\vspace{0.5cm}

\section{\hspace{0.2cm}Sketch of Methods}\label{sec:methods}

\vspace{0.3cm}

\textit{All integrals recorded in the above table without exceptions can be reduced to  linear combinations of \[\displaystyle \int_{0}^{+\infty} 
          \frac{\sinh\!\left((2k+1)x\right)}{x \cosh^{n}\! x} \, dx\] with rational coefficients. The derivations are detailed in \textnormal{\cite{talla_waffo_integral_2025}} and omitted here for the sake of brevity. The following lines only give instructions on how to prove that \[\displaystyle \int_{0}^{+\infty} 
          \frac{\sinh\!\left((2k+1)x\right)}{x \cosh^{2n+1}\! x} \, dx\] is of the form \textnormal{\eqref{eq:formZeta}} and $\displaystyle \int_{0}^{+\infty} 
          \frac{\sinh\!\left((2k+1)x\right)}{x \cosh^{2n}\! x} \, dx$ is of the form \textnormal{\eqref{eq:formBeta}}}

\vspace{0.3cm}

 The study of the integrals $\displaystyle \int_{0}^{+\infty} 
          \frac{\sinh\!\left((2k+1)x\right)}{x \cosh^{2n}\! x} \, dx$ and $\displaystyle \int_{0}^{+\infty} 
          \frac{\sinh\!\left((2k+1)x\right)}{x \cosh^{2n+1}\! x} \, dx$ relies solely on contour integration techniques by bearing in mind first of all the evenness of the integrands which imply that the integration bounds can be extended to $\left(-\infty, +\infty\right)$. We made the remark that the Laurent expansion approach is very fruitful to find the residues here. To find these residues, we note that $\cosh^{n}\! z$ has a pole of $n$-th order at $z_l := \dfrac{(2l+1)i\pi}{2}$ where $l \in \mathbb{Z}$.
We find the Laurent expansion of $\dfrac{1}{\cosh^{n}\! z}{}$ around $z_l$ by letting first of all $z = z_l + w$. It comes
\[\cosh^{n}\! z = (\cosh z_l \, \cosh w + \sinh z_l \, \sinh w)^n = \sinh^n z_l \, \sinh^n w\] since $\cosh z_l = 0$. Hence, finding the Taylor expansion of $\cosh^{n}\! z$ around $z_l$ reduces to finding the Taylor expansion of $\sinh^n z_l \, \sinh^n w$ around $w=0$, which is done after linearizing $\sinh^n w$ depending of the parity of $n$ and using the Taylor expansion of $\sinh(\alpha w)$ around $w=0$. With the Cauchy's product formula we exhibit the Laurent expansion of $\cosh^{-n}\! z$ around $z_l$ and remark rational coefficients appearing in these series. Since $\dfrac{\sinh\!\big((2k+1)z\big)}{z}$ is analytic near $z=z_l$, its Taylor expansion follows from Leibniz’s rule for derivatives of products together with the following identities, which are readily proved by induction on $n$.

\[ \dfrac{d^n}{dx^n}(x^{-1}) = (-1)^n \, n!\, x^{-(n+1)}\]
\[
\frac{d^n}{dx^n} \sinh(ax) = a^n \times 
\begin{cases}
\cosh(ax) & \text{if } n \text{ is odd}, \\
\sinh(ax) & \text{if } n \text{ is even}.
\end{cases}
\] The residue is read after multiplying the Laurent expansion of $\cosh^{-n}\! z$ by the Taylor expansion of $\dfrac{\sinh\!\left((2k+1)z\right)}{z}$. The contour of integration, designated by $C_{N, r}$ is the rectangle bounding the region of the complex plane defined by $\begin{cases} 
0 \leq \, Im(z) \leq 2N\pi\,\, (N \in \mathbb{N}) \\
\left|Re\, (z)\right| \leq r\,\, (r \in \mathbb{R_+^*})
\end{cases}$. The modulus of the integrand vanishes on the paths $\begin{cases} 
0 \leq \, Im(z) \leq 2N\pi\,\, \\
\left|Re\, (z)\right| = r\,\, 
\end{cases}$ and on the axis $\begin{cases} 
 Im(z) = 2N\pi\,\, \\
\left|Re\, (z)\right| \leq r\,\, 
\end{cases}$ for sufficiently large $N$ and $r$. Cauchy's residue theorem confirms the claim at the end.

\vspace{0.5cm}

The proofs of Malmsten’s integrals outlined in \cref{prop:zeta_polylog} and \cref{prop:beta_polylog} lack full mathematical rigor for two main reasons. First, the interchange of summation and integration is not properly justified. Second, the manipulations involve regularized versions of the beta and zeta functions, whose corresponding real-valued expressions would otherwise diverge. To establish a rigorous proof of these formulae, it is therefore necessary to derive explicit expressions for the coefficients $L_{p,n}$ and $K_{p,n}$ satisfying the given relations :
\[\int_0^1 \frac{\mathrm{Li}_{-2n-1}(-x^2) \, \ln\ln \frac{1}{x}}{x} \, dx = \sum_{p=1}^{n} L_{p,n} \, \frac{\zeta(2p+1)}{\pi^{2p}} \quad \text{and } \quad \int_0^1\dfrac{\operatorname{Im}\!\bigl[\operatorname{Li}_{-2n}(i\,x)\bigr]}{x} \,\ln \ln \frac{1}{x} \,dx = \sum_{p=1}^{n} K_{p,n} \, \frac{\beta(2p)}{\pi^{2p-1}}\]

, to compute the coefficients $L_{n,n}$ and $K_{n,n}$ and to show that \[L_{p,n} = K_{p,n} = 0 \qquad \forall p \in \llbracket 1, n-1\rrbracket\]

This method is fully explored and detailed in \cite[86-93, 123-130]{talla_waffo_integral_2025}. It starts by showing first of all the following identities :
\[
\begin{cases}
\displaystyle 
\mathrm{Li}_{-n}(z)
= \frac{1}{(1-z)^{n+1}}
  \sum_{k=0}^{n-1} 
  \genfrac{\langle}{\rangle}{0pt}{}{n}{k}\, z^{\,n-k}, \\[1.2em]
\displaystyle 
P_n(x)
:= \frac{1}{(1-x)^{n+1}}
   \sum_{k=0}^{n} 
   \EulerB{n}{k}\, x^{k}
= \sum_{k=0}^{\infty} (2k+1)^{n} x^{k}, \\[1.2em]
\displaystyle 
P_{2n}\!\bigl(-x^{2}\bigr)
= \frac{\operatorname{Im}\!\bigl[\operatorname{Li}_{-2n}(i\,x)\bigr]}{x},
\qquad n\in\mathbb{N},\; x\in\mathbb{R},\; |x|<1.
\end{cases}
\]

where $\displaystyle\genfrac{\langle}{\rangle}{0pt}{}{n}{k}$ are Eulerian numbers of type A and $\displaystyle\EulerB{n}{k}$ are Eulerian numbers of type B.

\vspace{0.5cm}

\section{\hspace{0.2cm}Conjectures and Perspectives}\label{sec:conjectures}

\vspace{0.5cm}

By virtue of \cref{prop:zeta_polylog}, the irrationality of $\dfrac{\zeta(2n+1)}{\pi^{2n+1}}$ is \textit{\textbf{equivalent}} to this conjecture :

The set
\[
\left\{ 
\pi, \ 
\int_0^1 
\frac{\operatorname{Li}_{-2n-1}(-x^2)}{x} \,\ln \ln \frac{1}{x} \,dx \ 
\right\}
\]
is linearly independent over $\mathbb{Z}$ with $n \geq 1$

\hspace{0.3cm}

By virtue of \cref{prop:beta_polylog}, the irrationality of $\dfrac{\beta(2n)}{\pi^{2n}}$ is \textit{\textbf{equivalent}} to this conjecture :

The set
\[
\left\{ 
\pi, \ 
\int_0^1 
\frac{\operatorname{Im}\!\bigl[\operatorname{Li}_{-2n}(i\,x)\bigr]}{x} \,\ln \ln \frac{1}{x} \,dx \ 
\right\}
\]
is linearly independent over $\mathbb{Z}$ with $n \geq 1$

\textit{In further lines of this paragraph, we mean by $\widetilde{\beta\zeta}(n)$ either the number \(\dfrac{\zeta(2n+1)}{\pi^{2n}}\) or \(\dfrac{\beta(2n)}{\pi^{2n-1}}\) with $n \in \mathbb{N}^*$}

\vspace{0.5cm}

\begin{proposition}\label[proposition]{prop:matrix_reverse}
Let $\mathbb{K}$ be a field, $V$ a $\mathbb{K}$-vector space, $I_n, J_n$ two sequences in $V$ such that holds :
\[ \forall n \in \mathbb{N}^*, \exists (x_{1,n}, x_{2,n}, \ldots , x_{n,n}) \in \mathbb{K}^n : I_n = \sum_{k=1}^n x_{k,n} J_k. \] Then \[\forall n \in \mathbb{N}^*, \exists (y_{1,n}, y_{2,n}, \ldots , y_{n,n}) \in \mathbb{K}^n : J_n = \sum_{k=1}^n y_{k,n} I_k \]
\end{proposition}

\begin{proof}
The exercise consists of expressing each $J_k$ as linear combination of $I_k$ which corresponds to the last step of Gauss elimination algorithm. The result is another formulation of the following statement :"\textit{The inverse of a lower triangular matrix is again a lower triangular matrix.}"\end{proof}

\vspace{0.7cm}

\textbf{\underline{Corollary :}} \textit{Since there exist integrals $I_n$ either holding form \textnormal{\eqref{eq:formZeta} or \eqref{eq:formBeta}}, we can always find rational coefficients $y_{k,n}$ satisfying $\widetilde{\beta\zeta}(n) = \displaystyle\sum_{k=1}^n y_{k,n} I_k$. The irrationality of $\dfrac{\widetilde{\beta\zeta}(n)}{\pi}$ would follow then from this conjecture :
\[
\pi \notin \text{Span}_{\mathbb{Q}}\left\{
I_1,\ 
I_2,\ 
I_3,\ 
I_4,\ 
I_5,\ 
I_6,\ 
I_7,\ 
I_8,\ 
I_9,\ 
I_{10},\ 
I_{11},\ 
I_{12},\ 
\ldots \right\}
\]
}
\vspace{0.3cm}

\textbf{\underline{Application to conjectures :}}

\vspace{0.3cm}

\begin{itemize}
	\item From the table of \cref{sec:result} and in combination with \cref{prop:matrix_reverse}, we conclude the existence of rational coefficients $\chi_{p,n}$ and $\xi_{p,n}$ which satisfy
\[
\dfrac{\zeta(2n+1)}{\pi^{2n}} = \displaystyle\sum_{p=1}^{n} \chi_{p,n} \displaystyle\int_{0}^{+\infty} \frac{\sinh x}{x\,\cosh^{2p+1}x}\,dx \quad \text{and } \quad
\dfrac{\beta(2n)}{\pi^{2n-1}} = \displaystyle\sum_{p=1}^{n} \xi_{p,n} \displaystyle\int_{0}^{+\infty} \frac{\sinh x}{x\,\cosh^{2p}x}\,dx\] The irrationality of the numbers \(\dfrac{\zeta(2n+1)}{\pi^{2n+1}}\) and \(\dfrac{\beta(2n)}{\pi^{2n}}\) would follow from this \textit{single} conjecture:

\[
\pi \notin \text{Span}_{\mathbb{Q}}\left\{
\int_{0}^{+\infty} \frac{\sinh x}{x\,\cosh^{2}x}\,dx,\ 
\int_{0}^{+\infty} \frac{\sinh x}{x\,\cosh^{3}x}\,dx,\ 
\int_{0}^{+\infty} \frac{\sinh x}{x\,\cosh^{4}x}\,dx,\ 
\int_{0}^{+\infty} \frac{\sinh x}{x\,\cosh^{5}x}\,dx,\ 
\ldots \right\}
\]

\item Putting the focus on the integral $\displaystyle\int_{0}^{1} \frac{x^{2p-1}}{\sqrt{1-x^2}\,\operatorname{arctanh }\, x}\,dx$, \cref{prop:matrix_reverse} ensures the existence of some rational coefficients $r_{p,n}$ such that
$
\dfrac{\beta(2n)}{\pi^{2n-1}} = \displaystyle\sum_{p=1}^{n} r_{p,n} \displaystyle\int_{0}^{1} \frac{x^{2p-1}}{\sqrt{1-x^2}\,\operatorname{arctanh }\, x}\,dx
$; thus, the existence of polynomials $\Xi_n \in \mathbb{Q}[x]$ such that $\dfrac{\beta(2n)}{\pi^{2n-1}} = \displaystyle\int_{0}^{1} \frac{x\, \Xi_n(x)}{\sqrt{1-x^2}\,\text{arctanh } x}\,dx$.

Few first polynomials $\Xi_n$ are \cite{talla_waffo_integral_2025}:

\vspace{0.5cm} 

\[
\Xi_{1}(x) = \frac{1}{4}
\]

\vspace{0.3cm}

\[
\Xi_{2}(x) = - \frac{1}{16}\,x^{2}
+ \frac{5}{96}
\]

\vspace{0.3cm}

\[
\Xi_{3}(x) = \frac{1}{64}\,x^{4}
- \frac{3}{128}\,x^{2}
+ \frac{61}{7 \, 680}
\]

\vspace{0.3cm}

\[
\Xi_{4}(x) = - \frac{1}{256}\,x^{6}
+ \frac{13}{1536}\,x^{4}
- \frac{173}{30\,720}\,x^{2}
+ \frac{277}{258 \, 048}
\]

\vspace{0.3cm}

\[
\Xi_{5}(x) = \frac{1}{1 \, 024}\,x^{8}
- \frac{17}{6 \, 144}\,x^{6}
+ \frac{203}{73 \, 728}\,x^{4}
- \frac{3\,403}{3\,096\,576}\,x^{2}
+ \frac{50\,521}{371\,589\,120}
\]

After multiplying each polynomial in $\mathbb{Q}[x]$ by the lowest common multiple of denominators of its coefficients, we get a polynomial in $\mathbb{Z}[x]$. Hence, the irrationality of the numbers \(\dfrac{\beta(2n)}{\pi^{2n}}\) for each $n$ would follow from the following conjecture:
\[\forall P \in \mathbb{Z}[x]\backslash \left\{0\right\}, \displaystyle\int_{0}^{1} \frac{x\, P(x)}{\sqrt{1-x^2}\,\text{arctanh } x}\,dx \notin \pi \mathbb{Q}\]

\item Putting the focus on the integral $\displaystyle\int_{0}^{1} \frac{x^{2p-1}}{\operatorname{arctanh }\, x}\,dx$, \cref{prop:matrix_reverse} ensures the existence of some rational coefficients $s_{p,n}$ such that
$
\dfrac{\zeta(2n+1)}{\pi^{2n}} = \displaystyle\sum_{p=1}^{n} s_{p,n} \displaystyle\int_{0}^{1} \frac{x^{2p-1}}{\operatorname{arctanh }\, x}\,dx
$; thus, the existence of polynomials $\Lambda_n \in \mathbb{Q}[x]$ such that $\dfrac{\zeta(2n+1)}{\pi^{2n}} = \displaystyle\int_{0}^{1} \frac{x\, \Lambda_n(x)}{\text{arctanh } x}\,dx$.

Few first polynomials $\Lambda_n$ are \cite[131-132]{talla_waffo_integral_2025}:

\vspace{0.5cm} 

\[
\Lambda_{1}(x) = \frac{1}{7}
\]

\vspace{0.3cm}

\[
\Lambda_{2}(x) = - \frac{1}{31}\,x^{2}
+ \frac{2}{93}
\]

\vspace{0.3cm}

\[
\Lambda_{3}(x) = \frac{1}{127}\,x^{4}
- \frac{4}{381}\,x^{2}
+ \frac{17}{5 \, 715}
\]

\vspace{0.3cm}

\[
\Lambda_{4}(x) = - \frac{1}{511}\,x^{6}
+ \frac{2}{511}\,x^{4}
- \frac{6}{2 \, 555}\,x^{2}
+ \frac{62}{160 \, 965}
\]

\vspace{0.3cm}

\[
\Lambda_{5}(x) = \frac{1}{2 \, 047}\,x^{8}
- \frac{8}{6 \, 141}\,x^{6}
+ \frac{37}{30 \, 705}\,x^{4}
- \frac{848}{1 \, 934 \, 415}\,x^{2}
+ \frac{1 \, 382}{29 \, 016 \, 225}
\]

After multiplying each polynomial in $\mathbb{Q}[x]$ by the lowest common multiple of denominators of its coefficients, we get a polynomial in $\mathbb{Z}[x]$. Hence, the irrationality of the numbers \(\dfrac{\zeta(2n+1)}{\pi^{2n+1}}\) for each $n$ would follow from the following conjecture:
\[\forall P \in \mathbb{Z}[x]\backslash \left\{0\right\}, \displaystyle\int_{0}^{1} \frac{x\, P(x)}{\text{arctanh } x}\,dx \notin \pi \mathbb{Q}\]
\end{itemize}

\textbf{\underline{In the interest of the theory of Fourier series : }}

\vspace{0.5cm}

Since \cref{prop:matrix_reverse} shows that $\displaystyle\dfrac{\zeta(2n+1)}{\pi^{2n}} = \sum_{p=1}^n r_{p,n} \int_0^{\frac{\pi}{4}} \frac{\sin(4px)}{\ln(\tan x)} \, dx$ and that \[\forall n \in \mathbb{N}^*, \dfrac{\beta(2n)}{\pi^{2n-1}} = \sum_{p=1}^n \Psi_{p,n} \int_0^{\frac{\pi}{4}} \frac{\cos((4p-2)x)}{\ln(\tan x)} \, dx\], with $r_{p,n}, \Psi_{p,n} \in \mathbb{Q}$, this conjecture, if true, would provide a pathway toward establishing the irrationality of $\dfrac{\beta(2n)}{\pi^{2n}}$ and $\dfrac{\zeta(2n+1)}{\pi^{2n+1}}$:

\[
\pi \notin \text{Span}_{\mathbb{Q}} \left\{
\int_0^{\frac{\pi}{4}} \frac{\cos(2x)}{\ln(\tan x)} \, dx, \ 
\int_0^{\frac{\pi}{4}} \frac{\sin(4x)}{\ln(\tan x)} \, dx, \ 
\int_0^{\frac{\pi}{4}} \frac{\cos(6x)}{\ln(\tan x)} \, dx, \ 
\int_0^{\frac{\pi}{4}} \frac{\sin(8x)}{\ln(\tan x)} \, dx, \
\ldots 
\right\}
\]

\vspace{0.5cm}

We establish an interesting identity, which lends positive support to the validity of this conjecture. \textit{Of course we recognize the weakness and lack of rigor of our lines of reasoning, because of the convergence and summation/integration swap.} We start from the well-known Fourier series expansion of $\ln \tan x$

\[
\ln(\tan x) \;=\; -2 \sum_{k=0}^{\infty} \frac{\cos\!\bigl((4k+2)x\bigr)}{2k+1},
\qquad 0 < x < \tfrac{\pi}{2}.
\]

Dividing both sides by $\ln \tan x$ -despite not rigorous justification of convergence, one has

\[
1 \;=\; -2 \sum_{k=0}^{\infty} 
\frac{\cos\!\bigl((4k+2)x\bigr)}{(2k+1)\,\ln(\tan x)},
\qquad 0 < x < \tfrac{\pi}{2}.
\]

On integrating term by term -despite lack of rigor, we obtain the formal identity

\[
\frac{\pi}{4}
\;=\;
-2 \sum_{k=0}^{\infty} \frac{1}{2k+1} 
\int_{0}^{\pi/4} \frac{\cos\!\bigl((4k+2)x\bigr)}{\ln(\tan x)} \, dx.
\]

Some computations of the first few partial sums on the right-hand side series tend to show that this identity may be correct. This identity is not consistent with the linear dependence of the previous set stated in our last conjecture.

\vspace{0.5cm}

\printbibliography

\end{document}